\documentclass[reqno]{amsart}
\usepackage{amscd,amsfonts,amssymb}
\usepackage{graphicx}
\usepackage{color}

\numberwithin{equation}{section}

\newtheorem{Theorem}{Theorem}[section]
\newtheorem{Lemma}{Lemma}[section]

\theoremstyle{definition}
\newtheorem{Definition}{Definition}[section]
\newtheorem{Proposition}{Proposition}[section]
\theoremstyle{remark}
\newtheorem{Remark}{Remark}[section]
\newtheorem{Example}{Example}[section]

\author{ Hailiang Liu  and Zhaoyang Yin}
\address{Iowa State University, Mathematics Department, Ames, IA 50011} \email{hliu@iastate.edu}
\address{Department of Mathematics, Sun Yat-sen University, 510275, Guangzhou,
China} \email{mcsyzy@mail.sysu.edu.cn}

\title[On nonlocal dispersive equations]{Global regularity, and wave breaking phenomena in a class of nonlocal dispersive equations}
\thanks{ }
\keywords{Dispersive equations, global regularity, wave breaking,
global weak solutions}
\date{\today}
\begin{document}
\begin{abstract}This paper is concerned with a class of nonlocal dispersive models -- the $\theta$-equation proposed by H. Liu [
On discreteness of the Hopf equation,  {\it Acta Math. Appl. Sin.} Engl. Ser. {\bf 24}(3)(2008)423--440]:
$$
(1-\partial_x^2)u_t+(1-\theta\partial_x^2)\left(\frac{u^2}{2}\right)_x
=(1-4\theta)\left(\frac{u_x^2}{2}\right)_x,
$$
including integrable equations such as the Camassa-Holm equation,
$\theta=1/3$, and the Degasperis-Procesi equation, $\theta=1/4$, as
special models. We investigate both global regularity of solutions
and wave breaking phenomena for $\theta  \in \mathbb{R}$. It is
shown that as $\theta$ increases regularity of solutions improves:
(i) $0 <\theta < 1/4$, the solution will blow up when the momentum
of initial data satisfies certain sign conditions; (ii) $1/4 \leq
\theta < 1/2$, the solution will blow up when the slope of initial
data is negative at one point; (iii) $\frac{1}{2} \leq \theta \leq
1$ and $\theta=\frac{2n}{2n-1}, n\in \mathbb{N}$, global existence
of strong solutions is ensured. Moreover, if the momentum of initial
data has a definite sign, then for any $\theta\in \mathbb{R}$ global
smoothness of the corresponding solution is proved. Proofs are
either based on the use of some global invariants or based on
exploration of favorable sign conditions of quantities involving
solution derivatives. Existence and uniqueness results of global
weak solutions for any $\theta \in \mathbb{R}$ are also presented.
For some restricted range of parameters results here are equivalent
to those known for the $b-$equations [e.g. J. Escher and Z. Yin,
Well-posedness, blow-up phenomena, and global solutions for the
b-equation, {\it J. reine angew. Math.}, {\bf 624} (2008)51--80.]
\end{abstract}
\maketitle \tableofcontents

\section{Introduction}
In recent years nonlocal dispersive models have been investigated
intensively at different levels of treatments: modeling, analysis as
well as numerical simulation.  The model derives in several ways, for instance,
(i) the asymptotic modeling of shallow water waves \cite{Wh, D-G-H1, D-G-H2}; (ii)
renormalization of dispersive operators \cite{Wh, Liu1}; and (iii) model equations
of some dispersive schemes \cite{Liu2}. The peculiar feature of
nonlocal dispersive models is their ability to capture both global
smoothness of solutions and the wave breaking phenomena.

In this work we focus on a class of nonlocal dispersive models --
the $\theta$-equation of the form
\begin{equation}\label{eq2.1}
(1-\partial_x^2)u_t+(1-\theta\partial_x^2)\left(\frac{u^2}{2}\right)_x
=(1-4\theta)\left(\frac{u_x^2}{2}\right)_x,
\end{equation}
subject to the initial condition
\begin{equation}\label{ini}
u(0, x)= u_0(x), \quad x\in \mathbb{R}.
\end{equation}
The equation can be formally rewritten as
\begin{equation}\label{eq2.1+}
u_t -u_{txx} +uu_x=\theta uu_{xxx}+(1-\theta)u_xu_{xx},
\end{equation}
which when $0<\theta<1$ involves a convex combination of nonlinear
terms $uu_{xxx}$ and $u_xu_{xx}$. This class was identified by H. Liu
\cite{Liu2} in his study of model equations for some dispersive
schemes to approximate the Hopf equation
$$
u_t+ uu_x=0.
$$
The model (\ref{eq2.1}) under a transformation links to the so
called $b-$model,
$$
u_t-\alpha^2 u_{txx} + c_0u_x +(b+1) u u_x +\Gamma u_{xxx}=\alpha^2
\left( b u_x u_{xx} + u u_{xxx}\right).
$$
which has been extensively studied in recent years
\cite{D-H-H1,D-P,E-Y, E-Y2,H-S,H-S2}. Both classes of equations
are contained in the more general class derived in \cite{Liu1}
using renormalization of dispersive operators and number of
conservation laws.

In (\ref{eq2.1}), two equations are worth of special attention:
$\theta=\frac{1}{3}$ and $\theta=\frac{1}{4}$. The $\theta$-equation
when $\theta=\frac{1}{3}$ reduces to the Camassa-Holm equation,
modeling the unidirectional propagation of shallow water waves over
a flat bottom, in which $u(t,x)$ denotes the fluid velocity at time
$t$ in the spatial $x$ direction \cite{C-H, D-G-H,J}. The
Camassa-Holm equation is also a model for the propagation of axially
symmetric waves in hyperelastic rods \cite{C-S2, Dai}. It has a
bi-Hamiltonian structure \cite{ F-F,L1} and is completely integrable
\cite{C-H, Cp}. Its solitary waves are smooth if $c_{0}>0$ and
peaked in the limiting case $c_{0}=0$, cf. \cite{C-H-H}. The orbital
stability of the peaked solitons is proved in \cite{C-S}, and that
of the smooth solitons in \cite{C-S3}. The explicit interaction of
the peaked solitons is given in \cite{B-S-S}.

The Cauchy problem for the Camassa-Holm equation has been studied
extensively. It has been shown that this problem is locally
well-posed \cite{C-Ep, Rb} for initial data $u_{0}\in
H^{3/2+}(\mathbb{R})$. Moreover, it has global strong solutions
\cite{Cf, C-Ep} and also admits finite time blow-up solutions
\cite{Cf, C-Ep, C-E}. On the other hand, it has global weak
solutions in $ H^1(\mathbb{R}) $ \cite{B-C, C-Ei,C-M,X-Z}. The
advantage of the Camassa-Holm equation in comparison with the KdV
equation,
$$
u_t +uu_x+\Gamma u_{xxx}=0,
$$
lies in the fact that the Camassa-Holm equation has peaked solitons
and models the peculiar wave breaking phenomena \cite{C-H-H,C-E}.

Taking $\theta=\frac{1}{4}$ in (\ref{eq2.1}) we find the
Degasperis-Procesi equation \cite{D-P}. The Degasperis-Procesi
equation can be regarded as a model for nonlinear shallow water
dynamics and its asymptotic accuracy is the same as that for the
Camassa-Holm shallow water equation \cite{D-G-H1,D-G-H2}. An inverse
scattering approach for computing $n$-peakon solutions to the
Degasperis-Procesi equation was presented in \cite{L-S}. Its
traveling wave solutions were investigated in \cite{L,V-P}. The
formal integrability of the Degasperis-Procesi equation was obtained
in \cite{D-H-H} by constructing a Lax pair. It has a bi-Hamiltonian
structure with an infinite sequence of conserved quantities and
admits exact peakon solutions which are analogous to the
Camassa-Holm peakons \cite{D-H-H}.

The study of the Cauchy problem for the Degasperis-Procesi equation
is more recent. Local well-posedness of this equation is established
in \cite{Y3} for initial data $u_{0}\in H^{3/2+}(\mathbb{R})$.
Global strong solutions are proved in \cite{E-L-Y,L-Y,Y4} and finite
time blow-up solutions in \cite{E-L-Y,L-Y,Y3,Y4}. On the other hand,
it has global weak solutions in $H^1(\mathbb{R})$, see e.g.
\cite{E-L-Y,Y4} and global entropy weak solutions belonging to the
class $L^{1}(\mathbb{R})\cap BV(\mathbb{R})$ and to the class
$L^{2}(\mathbb{R})\cap L^{4}(\mathbb{R})$, cf. \cite{Co-K}.

Though both the Dgasperis-Procesi and the Camassa-Holm equation share some nice properties, they differ in that the DP equation has not only peakon solutions \cite{D-H-H} and periodic peakon solutions \cite{Y5}, but also shock
peakons \cite{Lu} and the periodic shock waves \cite{E-L-Y2}.

The main quest of this paper is to see how regularity of solutions
changes in terms of the parameter $\theta$. With this in mind we
present a relative complete picture of solutions of problem (\ref{eq2.1})-(\ref{ini}) for different choices of $\theta$.

\begin{Theorem}\label{th1.1}{\rm [{\bf Global regularity}]} Let $u_0\in H^{3/2+}(\mathbb{R})$ and $m_0:=
(1-\partial_x^2)u_0$.
\begin{itemize}
\item[i)] For any $\theta \neq 0$, if in addition $u_0 \in
L^1(\mathbb{R})$, and $m_0$ has a definite sign ($m_0 \leq 0$ or
$m_0 \geq 0$ for all $x \in \mathbb{R}$), then the solution remains
smooth for all time. Moreover, for all $(t,x)\in
\mathbb{R}_{+}\times \mathbb{R}$, we have
\begin{itemize}
\item[(1)]  $m(t,x)u(t, x)\geq 0$ and
$$\|
m_{0}\|_{L^{1}(\mathbb{R})}=\| m(t, \cdot)\|_{L^{1}(\mathbb{R})}=\|
u(t,\cdot) \|_{L^{1}(\mathbb{R})}=\|
u_{0}\|_{L^{1}(\mathbb{R})}.$$
\item[(2)] $\|
u_{x}(t,\cdot)\|_{L^{\infty}(\mathbb{R})}\leq \|
u_0\|_{L^1(\mathbb{R})}$ and $$\|
u(t,\cdot)\|_{L^{\infty}(\mathbb{R})}\leq\| u(t,\cdot)\|_{1}\leq
\frac{\sqrt{2}}{2}e^{\frac{|\frac{1}{\theta}-3|t}{2}\|
u_{0}\|_{L^{1}(\mathbb{R})}} \| u_{0}\|_{1}.
$$
\end{itemize}
\item[ii)] For $\frac{1}{2} \leq \theta \leq 1$, if in addition $u_0 \in W^{2,
\frac{\theta}{1-\theta}}(\mathbb{R})$, then the solution remains
smooth for all time.
\item[iii)] For $\theta=\frac{2n}{2n-1}\in(1, 2), n\in \mathbb{N}$, if in addition $u_0 \in W^{3, \frac{\theta}{\theta-1}}(R)$, then the solution remains smooth for all time.
\end{itemize}
\end{Theorem}
\begin{Remark}
The result stated in i) recovers the global existence result of strong
solutions to the Camassa-Holm equation in \cite{C-Ep} and the
Degasperis-Procesi equation in \cite{Y4}.
\end{Remark}
\begin{Theorem}\label{th1.2}{\rm [{\bf Blow up criterion}]} Let $u_0\in H^{3/2+}(\mathbb{R})$ and $m_0:=
(1-\partial_x^2)u_0$.
\begin{itemize}
\item[i)] For $0< \theta \leq \frac{1}{4}$ and a fixed $x^*$, if $u_0(x^*+x)=-u_0(x^*-x)$  and $(x-x^*) m_0(x) \leq 0$ for any $x\in \mathbb{R}$, then the solution must blow up in finite time strictly before $T^*= -\frac{1}{u_{x}(0, x^*)}$ provided $u_x(0, x^*)<0$.
\item[ii)] For $\frac{1}{4}\leq  \theta <\frac{1}{2}$, if $u_0(x^*+x)=-u_0(x^*-x)$  for any $x\in \mathbb{R}$  and $u_{x}(0, x^*)<0$ , then the solution
must blow up in finite time strictly before
$T^*=\frac{2\theta}{(2\theta -1)u_{x}(0, x^*)}$.
\end{itemize}
\end{Theorem}
\begin{Remark}
The result stated in Theorem \ref{th1.2} shows that strong
solutions to the $\theta$-equation (\ref{eq2.1})-(\ref{ini}) for
$0<\theta<\frac{1}{2}$ may blow up in finite time, while Theorem
\ref{th1.1} shows that in the case $\frac{1}{2}\leq \theta\leq 1$
every strong solution to the $\theta$-equation (\ref{eq2.1})
exists globally in time. This presents a clear picture for global
regularity and blow-up phenomena of solutions to the
$\theta$-equation for all $0<\theta \leq 1$.
\end{Remark}

We shall present main ideas for proofs of the above results, a further refined analysis could be done following
those presented in \cite{E-Y} for the $b-$equation.
Note that the case $\theta=0$ is a borderline case  and not covered
by the above results, we shall present a detailed account for this
case. For completeness, we  also present global existence results for weak
solutions to characterize peakon solutions to (\ref{eq2.1}) for any
$\theta\in \mathbb{R}$.

The rest of this paper is organized as follows. In \S 2, we
present some preliminaries including how the $\theta-$equation
relates to other class of dispersive equations, the local
well-posedness and some key quantities to be used in subsequent
analysis. In \S 3, we show how global existence of smooth
solutions is established. The ideas for deriving the precise
blow-up scenario is given in \S 4. A detailed account for the case
of $\theta=0$ is presented in \S 5. Two existence and uniqueness
results on global weak solutions and one example for peakon
solutions to (1.1) for any $\theta\in \mathbb{R}$ are given in \S
6.

\section{Preliminaries}
\subsection{The $\theta-$equation and its variants}
The $\theta$-equation of the form
\begin{equation}\label{cta}
    (1-\partial_x^2)u_t+(1-\theta\partial_x^2)\left(\frac{u^2}{2}\right)_x
=(1-4\theta)\left(\frac{u_x^2}{2}\right)_x,
\end{equation}
up to a scaling of $t \to \frac{t}{\theta}$ for $\theta\not=0$, can be rewritten into a class of  B-equations
\begin{equation}\label{B}
u_t+ uu_x + [Q*B(u, u_x)]_x=0,
\end{equation}
where $Q=\frac{1}{2}e^{-|x|}$ and
$$
B=\left(\frac{1}{\theta}-1\right)\frac{u^2}{2} +\left( 4- \frac{1}{\theta} \right)\frac{u_x^2}{2}.
$$
The $B-$class with $B$ being quadratic in $u$ and $u_x$ was derived
in \cite{Liu1} by using a renormalization technique and examining
number of conservation laws. In this $B-$class the Camasa-Holm
equation corresponds to $B(u, u_x) = u^2+u_x^2/2$; and the
Degasperis-Procesi equation corresponds to $B(u, u_x) = 3u^2/2$. The
local well-posedness for (\ref{B}) with initial data $u_0(x)$ was
established in \cite{Liu1}.

\begin{Theorem}\cite{Liu1} Suppose that $u_0 \in H^{3/2+}_x$ and $B(u, p)$ are quadratic functions in its arguments,
then there exists a time $T$ and a unique solution $u$ of (\ref{B})
in the space $C([0, T);H^{3/2+}(\mathbb{R}))\cap C^1([0, T);
H^{1/2+}(\mathbb{R}))$ such that $lim_{t\downarrow 0}u(t,
\cdot)=u_0(\cdot)$. If $T< \infty$ is the maximal existence time,
then
$$
lim_{ t\to T}{\rm sup}_{0 \leq \tau \leq t}\|u_x(\cdot,
\tau)\|_{L^\infty}(\Omega)=\infty,$$ where $\Omega=\mathbb{R}$ for
initial data decaying at far fields, or $\Omega=[0, \pi]$ for
periodic data.
\end{Theorem}
Wave breaking criteria are identified separately for several
particular models in class (\ref{B}), using their special
features, see \cite{Liu1} for further details.

For $\theta\not=0$, the class of $\theta$-equations can also be
transformed into the $b-$equation of the form
\begin{equation}\label{b}
u_t-\alpha^2 u_{txx} + c_0u_x +(b+1) u u_x +\Gamma u_{xxx}=\alpha^2
\left( b u_x u_{xx} + u u_{xxx}\right).
\end{equation}
In fact, if we set
\begin{align*}
u(t, x)&=c_0 \theta + \tilde u(\tau, z), \\
 z &=\alpha\left( x -\theta
\left(1+\frac{\Gamma}{\alpha^2}\right)t\right),  \\
\tau &=\alpha \theta t,
\end{align*}
then a straightforward calculation leads to
$$
(1-\alpha^2 \partial_z^2) \tilde u_\tau +c_0\tilde u_z
+\frac{1}{\theta} \tilde u \tilde u_z +\Gamma \tilde
u_{zzz}=\alpha^2 \left( \left( \frac{1}{\theta}-1\right)\tilde u_z
\tilde u_{zz} +\tilde u \tilde u_{zzz}\right).
$$
Setting
$$
\theta =\frac{1}{b+1}
$$
and changing variables $(\tilde u, \tau, z)$ back to $(u, t, x)$ we
thus obtain the so-called $b-$equation (\ref{b}).

Note that the $\theta-$equation does not include  the case $b=-1$,
which has been known un-physical. Also $\theta=0$ case is not in the
class of $B-$equation (\ref{B}) either.

\subsection{Local well-posedness and a priori estimates }
In order to prove our main results for different cases, we need to establish the
following local existence result.
\begin{Theorem}\label{th2.2}{\rm [{\bf Local existence}]} Let $u_0\in H^{3/2+}(\mathbb{R})$, then exists a $T=T(\theta, \|u_0\|_{3/2+})>0$
and a unique solution in
$$
C([0, T);H^{3/2+}(\mathbb{R}))\cap C^1([0, T); H^{1/2+}(\mathbb{R})).
$$
The solution depends continuously on the initial data, i.e. the
mapping
$$u_{0} \rightarrow u(\cdot,u_{0}):H^{s}(\mathbb{R}) \rightarrow
C([0,T);H^{s}(\mathbb{R}))\cap C^{1}([0,T);H^{s-1}(\mathbb{R})),
\quad s>3/2 $$ is continuous. Moreover, if $T<\infty$ then
$\lim_{t \rightarrow T}\|u(t,\cdot)\|_s=\infty$.
\end{Theorem}
The proof for $\theta \not=0$ follows from that for the $b-$
equation in \cite{E-Y}  or for the $B-$ equation in \cite{Liu2}.

Furthermore we have the following result.
\begin{Theorem}\label{th2.3}
Let $u_{0} \in H^{3/2+}(\mathbb{R})$ be given and assume that $T$ is
the maximal existence time of the corresponding solution to
(\ref{eq2.1}) with the initial data $u_{0}$. If there exists an
$M>0$ such that
$$
\| u_{x}(t,x)\|_{L^{\infty}(\mathbb{R})}\leq M,\quad t\in[0,T),
$$
then the $H^{s}(\mathbb{R})-$ norm of $u(t,\cdot)$ does not blow up
for  $t\in [0,T)$.
\end{Theorem}
 Let $u$ be the solution in $C([0,T);H^{s}(\mathbb{R}))\cap C^{1}([0,T);H^{s-1}(\mathbb{R}))$,
 it suffices to verify how $\|u(t, \cdot)\|_s$
 depends on $\|u_x(t, \cdot)\|_\infty$.  Here we could carry out a careful energy estimate to obtain a differential inequality of the form
$$
\frac{d}{dt}\|u(t, \cdot)\|_s \leq C \|u_x(t, \cdot)\|_\infty
\|u(t, \cdot)\|_s.
$$
The claim then follows from the Gronwall inequality. A detailed
illustration of such a procedure for the case $\theta=0$ will be
given in \S 5.
\begin{Remark}
This result is fundamental for us to prove or disprove the global
existence of strong solutions. More precisely, global existence
follows from a priori estimate on $\|u_x(t, \cdot)\|_\infty$, and
the finite time blow up of $\|u_x(t, \cdot)\|_\infty$ under certain
initial conditions reveals the wave breaking phenomena.
\end{Remark}
\section{Global regularity}
\subsection{Key invariants and favorable sign conditions} Let $T$ be the life span of the strong solution
$u \in C([0, T);H^{3/2+}(\mathbb{R}))\cap C^1([0, T); H^{1/2+}(\mathbb{R}))$. We now look at some key estimates valid for $t\in [0, T)$.
First since the $\theta$-equation is in conservative form, so
\begin{equation}\label{uc}
  \int_{\mathbb{R}} udx=\int_{\mathbb{R}} u_0dx.
\end{equation}
Let $m=(1-\partial_x^2)u$, then
\begin{equation}\label{um}
    u=(1-\partial_x^2)^{-1}m=Q*m,
\end{equation}
which implies
\begin{equation}\label{1c}
    \int_{\mathbb{R}} mdx = \int_{\mathbb{R}} udx =\int_{\mathbb{R}} u_0dx=\int_{\mathbb{R}} m_0 dx.
\end{equation}
Moreover the equation (\ref{eq2.1}) can be reformulated as
\begin{equation}\label{meq}
    m_t +\theta u m_x +(1-\theta)mu_x=0.
\end{equation}
For any $\alpha \in \mathbb{R}$, let $x=x(t, \alpha)$ be the curve
determined by
$$
\frac{d}{dt} x=\theta u(t, x), \quad x(0, \alpha)=\alpha
$$
for $t\in [0, T)$.  Then $F=\frac{\partial x}{\partial \alpha}$ solves
$$
\frac{d}{dt} F = \theta u_x F
$$
as long as $u$ remains a strong solution. Along the curve $x=x(t,
\alpha)$ we also have
$$
\frac{d}{dt}m = (\theta-1) u_x m.
$$
These together when canceling the common factor $u_x$ leads to the following global invariant:
\begin{equation}\label{lg}
    m(t, x(t, \alpha))F^{\frac{1}{\theta}-1}=m_0(\alpha), \quad \forall \alpha \in \mathbb{R}.
\end{equation}
From this Lagrangian identity we see that $m$ has a definite sign
once $m_0$ has. Correspondingly it follows from (\ref{um}) that $u$ has a definite sign
\begin{equation}\label{sign}
{\rm sign}(m)={\rm sign}(m_0)={\rm sign}(u)
\end{equation}
provided that $m_0$ has a definite sign on $\mathbb{R}$.

From (\ref{lg})  it follows
$$
\int_{\mathbb{R}} |m|^{\frac{\theta}{1-\theta}}(t, x(t, \alpha)) F d\alpha
=\int_{\mathbb{R}} |m|^{\frac{\theta}{1-\theta}}dx=\int_{\mathbb{R}} |m_0|^{\frac{\theta}{1-\theta}}dx,
$$
which yields the following estimate:
\begin{equation}\label{2c}
\frac{d}{dt}\int_{\mathbb{R}} |m|^{\frac{\theta}{1-\theta}}dx=0.
\end{equation}
Inspired by \cite{D-H-H1}  we identify  another conservation laws as follows
\begin{equation}\label{3c}
     \frac{d}{dt} \int_{\mathbb{R}} \left( (1-\theta)^2 m^{\frac{2-\theta}{\theta -1}}m_x^2 + \theta^2 m^{\frac{\theta}{\theta -1}}
    \right)dx=0.
\end{equation}
This conserved quantity will be used for some cases in the range $\theta>1$.

Multiplying (\ref{meq}) by $m=u-u_{xx}$, and integrating by
parts, we obtain
\begin{equation}\label{mux}
\frac{d}{dt}\int_{\mathbb{R}} m^2 dx=(3\theta -2)\int_{\mathbb{R}}
u_xm^2dx,
\end{equation}
which suggests that $\theta=\frac{2}{3}$ is a critical point for
the blow-up scenario. Note that
\begin{equation}\label{mu2}
\|u(t, \cdot)\|_2 \leq  \|m(t, \cdot)\|_{L^2} \leq \sqrt{2} \|u(t,
\cdot)\|_2.
\end{equation}
Both (\ref{mux}) and (\ref{mu2}) together enable us to conclude the
following
\begin{Theorem}\label{th3.1-}
Assume $u_{0} \in H^{3/2+}(\mathbb{R})$.  If $\theta=\frac{2}{3}$,
then every solution to (\ref{eq2.1})-(\ref{ini}) remains regular
globally in time. If $\theta <\frac{2}{3}$, then the solution will
blow up in finite time if and only if the slope of the solution
becomes unbounded from below in finite time. If $\theta
>\frac{2}{3}$, then the solution will blow up in finite time if and
only if the slope of the solution becomes unbounded from above in
finite time.
\end{Theorem}
\begin{Remark}
This result not only covers the corresponding results for the
Camassa-Holm equation in \cite{Cf,Y1} and the Degasperis-Procesi
equation in \cite{Y3}, but also presents another different possible
blow-up mechanism, i.e., if $\theta>\frac{2}{3}$, then the solution
to (\ref{eq2.1}) blows up in finite time if and only if the slope of
the solution becomes unbounded from above in finite time.
\end{Remark}

\subsection{Global existence: proof of Theorem \ref{th1.1} }  Let $T$ be the
maximum existence time of the solution $u$ with initial data
$u_0\in H^s$. Using a simple density argument we can just consider
the case $s=3$. Based on Theorem \ref{th2.2} and Theorem
\ref{th2.3} it suffices to show the uniform bound of $\|u_x(t,
\cdot)\|_\infty$ for all cases presented in Theorem \ref{th1.1}.

The proof of the first assertion i) is based on the global invariant (\ref{lg}), which implies (\ref{sign}), i.e.,  $m$ has a definite sign for $t>0$ as long as $m_0$ has a definite sign. Then for any $(t, x)\in [0, T)\times
\mathbb{R}$,
\begin{align*}
|u_x(t, x)|&=|Q_x*m| \leq \|Q_x\|_\infty
\|m\|_{L^1}=\frac{1}{2}\left|\int_{\mathbb{R}} m dx \right| \\
& =\frac{1}{2}\left|\int_{\mathbb{R}} u dx \right| =
\frac{1}{2}\left|\int_{\mathbb{R}} u_0 dx \right|\leq
\frac{1}{2}\|u_0\|_{L^1}.
\end{align*}
Then $T=\infty.$

The second assertion ii) follows from the use of (\ref{2c}), i.e.,
$$
\int_{\mathbb{R}}
|m|^{\frac{\theta}{1-\theta}}dx=\int_{\mathbb{R}}
|m_0|^{\frac{\theta}{1-\theta}}dx \leq  \|u_0\|_{W^{2, p}}, \quad
p=\frac{\theta}{1-\theta} \in [1, \infty].
$$
From $m\in L^p(\mathbb{R})$ and $u-u_{xx}=m$ it follows that $u\in
W^{2, p}(\mathbb{R})$. By the Sobolev imbedding theorem, we see that
$W^{2, p}(\mathbb{R}) \subset C^1(\mathbb{R})$. Thus $T=\infty$.

The last assertion iii) follows from the use of (\ref{3c}) with $\theta=\frac{2n}{2n-1}$, which upon integration leads to
$$
\int_{\mathbb{R}}( m^{2n-2}m_x^2+4n^2 m^{2n})dx=\int_{\mathbb{R}}(
m_0^{2n-2}m_{0x}^2+4n^2 m_0^{2n})dx.
$$
From this we see that  $m\in L^\infty$, for
$$
m^{2n}=\int_{-\infty}^x 2nm^{2n-1}m_x dx \leq \frac{1}{2}
\int_{\mathbb{R}}( m^{2n-2}m_x^2+4n^2 m^{2n})dx.
$$
Using $u=Q*m$ we obtain that $u\in W^{2, \infty}$; that is $ |u_x|$
is uniformly bounded. Thus $T=\infty$.
\section{Blow up phenomena: proof of Theorem \ref{th1.2}}
For the blow up analysis, one needs to find a way to show that
$d=u_x$ will become unbounded in  finite time. Rewriting (\ref{cta})  as
$$
u_t+\theta uu_x =\frac{Q_x}{2}*\left[(1-4\theta)u_x^2
+(\theta-1)u^2 \right].
$$
Notice that $Q_{xx}=Q-\delta(x)$;  a direct differentiation in $x$ of the
above equation leads to
$$
d_t+\theta ud_x +\left( \frac{1}{2} -\theta \right)d^2 =
\frac{1-\theta}{2}u^2+ \frac{Q}{2}*\left[ (1-4\theta)d^2
+(\theta-1)u^2 \right].
$$
For $\theta<\frac{1}{2}$ there is no control on $u^2$ term while we
track dynamics of $d$. The idea here, motivated by that used in
\cite{E-Y}, is to focus on a curve $x=h(t)$ such that $u(t, h(t))=0$
and $h(0)=x^*$.  On this curve
\begin{equation}\label{dh}
 \dot d +\left( \frac{1}{2} -\theta \right)d^2=\frac{Q}{2}*\left[
(1-4\theta)u_x^2 +(\theta-1)u^2 \right](t, h(t)).
\end{equation}
Two cases are distinguished: \\
(i) $\frac{1}{4} \leq \theta <\frac{1}{2}$. In this range of
$\theta$, the right-hand side of (\ref{dh}) is non-positive. We thus have
$$
\dot d +\left( \frac{1}{2} -\theta \right)d^2 \leq 0,
$$
for which $d$ will become unbounded from below in  finite time as
long as
$d(0, h(0))=u_x(0, x^*)<0$.\\
(ii) $0< \theta < \frac{1}{4}$. In this range of $\theta$ we also
need control the nonlocal term. If we can identify some initial data
such that
\begin{equation}\label{uux}
 Q* \left[ (1-4\theta)u_x^2 +(\theta-1)u^2\right](t, h(t)) \leq (1-4\theta) \left[u_x^2
-u^2 \right](t, h(t)) =(1-4\theta)d(t)^2.
\end{equation}
Then we have
$$
\dot d +\left( \frac{1}{2} -\theta \right)d^2 \leq \left(\frac{1}{2}-2\theta\right)d^2.
$$
That is
$$
\dot d +\theta d^2 \leq 0.
$$
Again in this case $d$ will become unbounded from below in finite
time once $d(0, h(0))=u_x(0, x^*)<0$.

Now we verify that the assumptions in Theorem \ref{th1.2} are
sufficient for claim (\ref{uux}) to hold. From
$u_0(x^*+x)=-u_0(x^*-x)$ for any $x\in \mathbb{R}$, it follows that
$u_0(x^*)=0$, and $u(t, x^*+x)=-u(t, x^*-x)$ due to symmetry of the equation. We then have $u(t,
x^*)=0$, leading to the case $h(t)=x^*$.

We further assume that
$$
(x-x^*)m_0(x) \leq 0,
$$
which combined with  (\ref{lg}) yields
$$
(x-x^*)m(t, x) \leq 0.
$$
This relation enables one  to use a similar argument as  (5.3)-(5.10) in \cite{E-Y}  to obtain
$$
Q*[u_x^2-u^2](t, x^*) \leq (u_x^2-u^2)(t, x^*).
$$
Hence
\begin{align*}
Q*\left[ (1-4\theta)u_x^2 +(\theta-1)u^2\right](t, x^*) & =(1-4\theta)Q*[u_x^2-u^2](t, x^*)-3\theta Q*[u^2](t, x^*)\\
& \leq (1-4\theta)Q*[u_x^2-u^2](t, x^*)\\
& \leq (1-4\theta)(u_x^2-u^2)(t, x^*),
\end{align*}
which leads to (\ref{uux}) as desired.

\section{A detailed account of the case $\theta=0$} In this
section, we establish the local well-posedness and present the
precise blow-up scenario and global existence results for the
$\theta-$equation with $\theta=0$, i.e.,
\begin{equation}\label{c0}
u_t -u_{txx}=u_xu_{xx}-uu_x.
\end{equation}
Note that  $(1- \partial^{2}_{x})^{-1}f = Q*f $ for all $f \in
L^{2}(\mathbb{R})$ and $Q \ast m=u $ for $m=u-u_{xx}$. Using this
relation, we can rewrite (\ref{c0}) as follows:
\begin{equation}\label{c0+}
\left\{\begin{array}{ll} u_{t}=\partial_{x}Q\ast
(\frac{1}{2}u_x^{2}-\frac{1}{2}u^{2}) , \quad &t
> 0,\;x\in \mathbb{R},\\ u(0,x) = u_{0}(x),&x\in \mathbb{R},\end{array}\right.
\end{equation}
or in the equivalent form:
\begin{equation}\label{c0++}
\left\{\begin{array}{ll} u_{t}=\partial_{x}(1-\partial^{2}_{x})^{-1}
(\frac{1}{2}u_x^{2}-\frac{1}{2}u^{2}),\quad &t>0,\; x\in \mathbb{R}, \\
u(0,x)=u_{0}(x),&x\in \mathbb{R}.
\end{array}\right.
\end{equation}
\begin{Theorem} \label{th3.1} Given $u_{0} \in H^{s}(\mathbb{R}),\;s>\frac{3}{2}$,
there exists a $T=T(\| u_{0}\|_{s})>0$, and a unique solution $u$
to (\ref{c0}) such that
$$
u=u(\cdot,u_{0})\in C([0,T);H^{s}(\mathbb{R}))\cap
C^{1}([0,T);H^{s}(\mathbb{R})).
$$
The solution depends continuously on the initial data, i.e. the
mapping $$u_{0} \rightarrow u(\cdot,u_{0}):H^{s}(\mathbb{R})
\rightarrow C([0,T);H^{s}(\mathbb{R}))\cap
C^{1}([0,T);H^{s}(\mathbb{R}))$$ is continuous. Moreover, if
$T<\infty$ then $\lim_{t \rightarrow T}\|u(t,\cdot)\|_s=\infty$.
\end{Theorem}
\begin{proof} Set $f(u)=\partial_{x}Q\ast
(\frac{1}{2}u_x^{2}-\frac{1}{2}u^{2})=\partial_{x}(1-\partial^2_{x})^{-1}
(\frac{1}{2}u_x^{2}-\frac{1}{2}u^{2})$. Let $u,v\in
H^{s},s>\frac{3}{2}$. Note that $H^{s-1}$ is a Banach algebra.
Then, we have
\begin{eqnarray}\label{3.4}
&\parallel&f(u)-f(v)
\parallel_{s}\nonumber\\&=&\parallel \partial_{x}(1-\partial^{2}_{x})^{-1}\left(\frac{1}{2}(u^{2}-v^{2})
+\frac{1}{2}(u^{2}_{x}-v^{2}_{x})\right)\parallel_{s}\nonumber\\
&\leq&\frac{1}{2}\parallel
(u-v)(u+v)\parallel_{s-1}+\frac{1}{2}\parallel
(u_{x}-v_{x})(u_{x}+v_{x})\parallel_{s-1}\\
&\leq&\frac{1}{2}\parallel u-v\parallel_{s}\parallel
u+v\parallel_{s}+\frac{1}{2}\parallel
\partial_{x}(u-v)\parallel_{s-1}\parallel u_{x}+v_{x}\parallel_{s-1}\nonumber\\
&\leq&\left(\parallel u
\parallel_{s}+\parallel v \parallel_{s}\right)\parallel u-v
\parallel_{s}.\nonumber
\end{eqnarray}
This implies that $f(u)$ satisfies a local Lipschitz condition in
$u$, uniformly in $t$ on $[0,\infty)$.

Next we show that for every $t_0\geq 0$, $u(t_0)\in
H^s(\mathbb{R})$, the Cauchy problem (\ref{c0}) has a unique mild
solution $u$ on an interval $[t_0, t_1]$ whose length is bounded
below by
$$
\delta(\|u(t_0)\|_s)=\frac{\|u(t_0)\|_s}{r^2(t_0)}=\frac{1}{2r(t_0)},
$$
where $r(t_0)=2\|u(t_0)\|_s$. Set $t_1=t_0+\delta(\|u(t_0)\|_s)$.
Let us define by
$$
\|u\|_{C([t_0,t_1];H^s(\mathbb{R}))}:=\sup_{t\in[t_0,t_1]}\|u\|_s
$$
the norm of $u$ as an element of $C([t_0,t_1];H^s(\mathbb{R}))$.
For a given $u(t_0)\in H^s$ we define a mapping
$F:C([t_0,t_1];H^s(\mathbb{R}))\longrightarrow
C([t_0,t_1];H^s(\mathbb{R}))$ by
\begin{equation}\label{3.5}
(Fu)(t)=u(t_0)+\int_{t_0}^t f(u(s))\,ds,\quad t_0\leq t\leq t_1.
\end{equation}
The mapping $F$ defined by (\ref{3.5}) maps the ball of radius $r(t_0)$
centered at $0$ of $C([t_0,t_1];H^s(\mathbb{R}))$ into itself.
This follows from the following estimate
\begin{equation}
\begin{split}
\|(Fu)(t)\|_s&\leq\|u(t_0)\|_s+\int_{t_0}^t \|f(u(s))-f(0)\|_s\,
ds
\\ &\leq\|u(t_0)\|_s+r^2(t_0)(t-t_0)
\\&\leq 2\|u(t_0)\|_s=r(t_0),
\end{split}
\end{equation}
where we have used the relations (\ref{3.4})-(\ref{3.5}), $f(0)=0$
and the definition of $t_1$.

By (\ref{3.4}) and (\ref{3.5}), we have
\begin{equation}\label{3.7}
\|(Fu)(t)-(Fv)(t)\|_s\leq
2r(t_0)(t-t_0)\|u-v\|_{C([t_0,t_1];H^s(\mathbb{R}))}.
\end{equation}
Using (\ref{3.5}) and (\ref{3.7}) and induction on $n$, we obtain
\begin{equation}\label{3.8}
\begin{split} \|(F^nu)(t)-(F^nv)(t)\|_s&\leq
\frac{\left(2r(t_0)(t-t_0)\right)^n}{n!}\|u-v\|_{C([t_0,t_1];H^s(\mathbb{R}))}\\&\leq
\frac{\left(2r(t_0)\delta(\|u(t_0)\|_s)\right)^n}{n!}\|u-v\|_{C([t_0,t_1];H^s(\mathbb{R}))}\\&\leq
\frac{1}{n!}\|u-v\|_{C([t_0,t_1];H^s(\mathbb{R}))}.
\end{split}\end{equation} For $n\geq 2$ we have $
\frac{1}{n!}<1$. Thus, by a well known extension of the Banach
contraction principle, we know that $F$ has a unique fixed point
$u$ in the ball of $C([t_0,t_1];H^s(\mathbb{R}))$. This fixed
point is the mild solution of the following integral equation
associated with Eq.(\ref{c0}):
\begin{equation}\label{3.9}
u(t,x)=u(t_0,x)+\int_{t_0}^t\partial_{x}Q\ast
(\frac{1}{2}u_x^{2}-\frac{1}{2}u^{2})(\tau,x)d\tau.
\end{equation}

Next, we prove the uniqueness of $u$ and the Lipschitz continuity
of the map $u(t_0)\longrightarrow u$. Let $v$ be a mild solution
to (\ref{c0}) on $[t_0,t_1)$ with initial data $v(t_0)$. Note that
$\|u\|_s\leq 2\|u(t_0)\|_s$ and $\|v\|_s\leq 2\|v(t_0)\|_s$. Then
\begin{equation}
\begin{split}
&\quad\,\,\|u(t)-v(t)\|_s\\&\leq
\|u(t_0)-v(t_0)\|_s+\int_{t_0}^{t}\|f(u)-f(v)\|_sd\tau
\\&\leq \|u(t_0)-v(t_0)\|_s+ (\|u(t_0)\|_s+\|v(t_0)\|_s)\int_{t_0}^t\|u(t)-v(t)\|_sd\tau.
\end{split}\end{equation} An application of Gronwall's inequality yields
\begin{equation*}
\|u(t)-v(t)\|_s\leq
e^{(\|u(t_0)\|_s+\|v(t_0)\|_s)(t_1-t_0)}\|u(t_0)-v(t_0)\|_s.
\end{equation*}
Therefore
\begin{equation}
\|u-v\|_{C([t_0,t_1];H^s(\mathbb{R}))}\leq
e^{(\|u(t_0)\|_s+\|v(t_0)\|_s)(t_1-t_0)}\|u(t_0)-v(t_0)\|_s,
\end{equation}
which implies both the uniqueness of $u$ and the Lipschitz
continuity of the map $u(t_0)\longrightarrow u$.

From the above we know that if $u$ is a mild solution of (\ref{c0}) on
the interval $[0,\tau]$, then it can be extended to the interval
$[0,\tau+\delta]$ with $\delta>0$ by defining on
$[\tau,\tau+\delta]$, $u(t,x)=v(t,x)$ where $v(t,x)$ is the
solution of the following integral equation
$$
v(t,x)=u(\tau)+\int_{\tau}^{t}f(v(s,x))ds,\qquad \tau\leq t\leq
\tau+\delta,
$$
where $\delta$ depends only on $\|u(\tau,\cdot)\|_s$. Let $T$ be
the maximal existence time of the mild solution $u$ of (\ref{c0}). If
$T<\infty$ then $\lim_{t\rightarrow T}\|u(t,\cdot\|_s)=\infty$.
Otherwise there is a sequence $t_n\longrightarrow T$ such that
$\|u(t_n,\cdot)\|_s\leq C$ for all $n$. This would yield that for
each $t_n$, near enough to $T$, $u$ defined on $[0,t_n]$ can be
extended to $[0,t_n+\delta]$ where $\delta>0$ is independent of
$t_n$. Thus $u$ can be extended beyond $T$. This contradicts the
definition of $T$.

Note that $u\in C([0,T);H^s(\mathbb{R}))$ and $f(u)$ satisfies
locally Lipschitz conditions in $u$, uniformly in $t$ on $[0,T)$.
Then we have that $f(u(t,x))$ is continuous in $t$. Thus it
follows from (\ref{3.7})-(\ref{3.9}) that
$$u(t,x)\in C([0,T);H^{s}(\mathbb{R}))\cap
C^{1}([0,T);H^{s}(\mathbb{R}))$$ is the solution to (\ref{c0}). This
completes the proof of the theorem.
\end{proof}

Next, we present the precise blow-up scenario for solutions to
Eq.(\ref{c0}).
\par
We first recall the following two useful lemmas.
\begin{Lemma}\label{lem3.1} \cite{K-P}
If $r>0$, then $H^{r}(\mathbb{R})\bigcap L^{\infty}(\mathbb{R})$
is an algebra. Moreover
$$
\| fg\|_{r}\leq c(\| f\|_{L^{\infty}(\mathbb{R})}\| g\|_{r}+\| f
\|_{r}\| g \|_{L^{\infty}(\mathbb{R})}),
$$
where $c$ is a constant depending only on $r$.
\end{Lemma}
\begin{Lemma}\label{lem3.2} \cite{K-P}
For $\Lambda=(1-\partial_x^2)^{1/2}$.  If $r>0$, then
$$
\| [\Lambda^{r},f]g\|_{L^{2}(\mathbb{R})}\leq c(\|
\partial_{x}f\|_{L^{\infty}(\mathbb{R})}\|
\Lambda^{r-1}g\|_{L^{2}(\mathbb{R})}+\| \Lambda^{r}f
\|_{L^{2}(\mathbb{R})}\| g \|_{L^{\infty}(\mathbb{R})}),
$$
where $c$ is a constant depending only on $r$.
\end{Lemma}
Then we prove the following useful result.
\begin{Theorem}\label{th3.2}
Let $u_{0} \in H^{s}(\mathbb{R}),\,s > \frac{3}{2}$ be given and
assume that $T$ is the existence time of the corresponding
solution to Eq.(\ref{c0}) with the initial data $u_{0}$. If there
exists $M>0$ such that
$$
\| u_{x}(t,x)\|_{L^{\infty}(\mathbb{R})}\leq M,\quad t\in[0,T),
$$
then the $H^{s}(\mathbb{R})-$ norm of $u(t,\cdot)$ does not blow
up on $[0,T)$.
\end{Theorem}
\begin{proof}  Let $u$ be the solution to Eq.(\ref{c0}) with initial data $u_{0}\in
H^{s}(\mathbb{R})$, $s>\frac{3}{2}$, and let $T$ be the maximal
existence time of the solution $u$, which is guaranteed by Theorem \ref{th3.1}.
Throughout this proof, $c>0$ stands for a generic constant
depending only on $s$.
\par
Applying the operator $\Lambda^{s}$ to Eq.(\ref{c0+}), multiplying by
$\Lambda^{s}u$, and integrating over $\mathbb{R}$, we obtain
\begin{equation}\label{3.12}
\frac{d}{dt}\| u\|_{s}^{2}=2(u,f_1(u))_{s}+2(u,f_2(u))_{s},
\end{equation}
where
$$f_1(u)=\partial_{x}(1-\partial^{2}_{x})^{-1}(-\frac{1}{2}u^{2})=-(1-\partial^{2}_{x})^{-1}(uu_x)$$
and $f_2(u)=\partial_{x}(1-\partial^{2}_{x})^{-1}
(\frac{1}{2}u^{2}_{x}).$
\\ Let us estimate the first term of the right-hand side of
Eq.(\ref{3.12}).
\begin{eqnarray}\label{eq3.3}
\mid (f_1(u),u)_{s}\mid &=&\mid
(\Lambda^{s}(1-\partial^{2}_{x})^{-1}(u\partial_{x}u),\Lambda^{s}u)_{0}\mid\nonumber
\\&\leq&\mid (\Lambda^{s-1}(u\partial_{x}u),\Lambda^{s-1}u)_{0}\mid
\nonumber\\&\leq&
\mid([\Lambda^{s-1},u]\partial_{x}u,\Lambda^{s-1}u)_{0}
+(u\Lambda^{s-1}\partial_{x}u,\Lambda^{s-1}u)_{0}\mid
\nonumber\\&\leq& \|
[\Lambda^{s-1},u]\partial_{x}u\|_{0}\|\Lambda^{s-1}u\|_{0}
+\frac{1}{2}\mid(u_{x}\Lambda^{s-1}u,\Lambda^{s-1}u)_{0}\mid
\nonumber\\&\leq& (c\| u_{x}
\|_{L^{\infty}(\mathbb{R})}+\frac{1}{2}\| u_{x}
\|_{L^{\infty}(\mathbb{R})})\| u \|^{2}_{s-1}\nonumber\\&\leq& c\|
u_{x} \|_{L^{\infty}(\mathbb{R})}\| u \|^{2}_{s}.
\end{eqnarray}
Here, we applied Lemma \ref{lem3.2} with $r=s-1$. Then, let us
estimate the second term of the right-hand side of (\ref{3.12}).
\begin{eqnarray}\label{eq3.4}
\mid (f_2(u),u)_{s}\mid&\leq&\| f_2(u)\|_{s}\| u
\|_{s}\leq\frac{1}{2}\| u^{2}_{x}\|_{s-1}\| u \|_{s}
\nonumber\\&\leq& c(\| u_{x} \|_{L^{\infty}(\mathbb{R})}\|
u\|_{s-1})\| u \|_{s} \nonumber\\&\leq& c\| u_{x}
\|_{L^{\infty}(\mathbb{R})}\| u \|^{2}_{s},
\end{eqnarray}
where we  used Lemma \ref{lem3.1} with $r=s$. Combining inequalities
(\ref{eq3.3})-(\ref{eq3.4}) with (\ref{3.12}), we obtain
\begin{equation*}
\frac{d}{dt}\| u \|^{2}_{s}\leq cM\| u\|^{2}_{s}.
\end{equation*}
An application of Gronwall's inequality yields
\begin{equation}
\| u(t) \|_{s}^{2} \leq\exp\left(cMt\right)\| u(0) \|_{s}^{2}.
\end{equation}
This completes the proof of the theorem.\end{proof}

We now present the precise blow-up scenario for Eq.(\ref{c0}).
\begin{Theorem} \label{th3.3}Assume that $u_{0} \in H^{s}(\mathbb{R}),\,s> \frac{3}{2}$.
Then the solution to Eq.(\ref{c0}) blows up in finite time if and only
if the slope of the solution becomes unbounded from below in
finite time.
\end{Theorem}
\par
\begin{proof} Applying Theorem \ref{th3.1} and a simple density argument, it suffices to consider the case $s=3$.
Let $ T > 0 $ be the maximal time of existence of the solution $u$
to Eq.(\ref{c0}) with initial data $u_{0}\in H^{3}(\mathbb{R})$.
From Theorem \ref{th3.1} we know that $u\in
C([0,T);H^{3}(\mathbb{R}))\cap C^{1}([0,T);H^{3}(\mathbb{R}))$.
\par
Multiplying Eq.(\ref{c0}) by $u$ and integrating by parts, we get
\begin {equation}\label{3.16}
\frac{d}{dt}\int_{\mathbb{R}}(u^{2}+u_x^2)dx=2\int_{\mathbb{R}}uu_{x}u_{xx}dx
-2\int_{\mathbb{R}}u^2u_{x}dx=-\int_{\mathbb{R}}u_{x}u_x^{2}dx.
\end{equation}
Differentiating Eq.(\ref{c0}) with respect to $x$, then multiplying the
obtained equation by $u_x$ and integrating by parts, we obtain
\begin {equation}\label{3.17}
\begin{split}
\frac{d}{dt}\int_{\mathbb{R}}(u_x^{2}+u_{xx}^2)dx&=-2\int_{\mathbb{R}}u_xu^2_{xx}dx
+2\int_{\mathbb{R}}uu_xu_{xx}dx\\&=-2\int_{\mathbb{R}}u_{x}u_x^{2}dx-\int_{\mathbb{R}}u_x^3dx.
\end{split}\end{equation}
Summing up (\ref{3.16}) and (\ref{3.17}), we have
\begin {equation}\label{3.18}
\frac{d}{dt}\int_{\mathbb{R}}(u^{2}+2u_x^2+u^2_{xx})dx=-\int_{\mathbb{R}}u_{x}(u_x^{2}+u_{xx}^2)dx.
\end{equation}
If the slope of the solution is bounded from below on $[0,T)\times
\mathbb{R}$, then there exists $M>0$ such that
\begin {equation*}
\frac{d}{dt}\|u\|_2\leq M\|u\|_2.
\end{equation*}
By means of Gronwall's inequality, we have
\begin{equation*}
\| u(t,\cdot)\|_{2}\leq  \| u(0,\cdot)\|_{2}\exp\{Mt\}, \quad
\forall t\in [0,T).
\end{equation*}
By Theorem \ref{th3.2}, we see that the solution does not blow up in
finite time.

On the other hand, by Theorem \ref{th3.1} and Sobolev's imbedding theorem,
we see that if the slope of the solution becomes unbounded from
below in finite time, then the solution will blow up in finite
time. This completes the proof of the theorem.\end{proof}
\begin{Remark}
Theorem \ref{th3.3} shows that (\ref{c0}) has the same blow-up scenario as
the Camassa-Holm equation \cite{Cf,Y1} and the Degasperis-Procesi
equation \cite{Y3} do.
\end{Remark}

Finally, we show that there exist global strong solutions to
Eq.(\ref{c0}) provided the initial data $u_0$ satisfies certain sign
conditions.

\begin{Lemma}\label{lem3.3}
Assume that $u_{0}\in H^{s}(\mathbb{R}),\;s> \frac{3}{2}$. Let
$T>0$ be the existence time of the corresponding solution $u$ to
(\ref{c0}). Then we have
\begin{equation}
m(t,x)=m_0(x)\exp^{-\int_0^tu_x(\tau,\,x)\,d\tau},
\end{equation}
where $(t,x)\in[0,T)\times \mathbb{R}$ and $m=u-u_{xx}$. Moreover,
for every $(t,x)\in[0,T)\times \mathbb{R}$, $m(t,x)$ has the same
sign as $m_0(x)$ does.
\end{Lemma}
\par
\begin{proof} Let $T>0$ be the maximal existence time of the solution $u$ with
initial data $u_0\in H^{s}(\mathbb{R})$.

Due to $u(t,x)\in C^{1}([0,T);H^{s}(\mathbb{R}))$ and
$H^{s}(\mathbb{R})\subset C(\mathbb{R})$, we see that the function
$u_{x}(t,x)$ are bounded, Lipschitz in the space variable $x$, and
of class $C^{1}$ in time. For arbitrarily fixed $T'\in(0,T)$,
Sobolev's imbedding theorem implies that
$$
\sup_{(s,x)\in[0,T']\times \mathbb{R}}\mid u_{x}(s,x)\mid<\infty.
$$
Thus, we infer from the above inequality that there exists a
constant $K>0$ such that
\begin{equation} \label{eq3.7}
e^{-\int_0^tu_x(\tau,\,x)\,d\tau}\geq e^{-tK}>0 \quad\text{ for
}\quad (t,x)\in[0,T']\times \mathbb{R}.
\end{equation}
By Eq.(\ref{c0}) and $m=u-u_{xx}$, we have
\begin{equation}
m_t(t,x)=-u_x(t,x)m(t,x).
\end{equation}
This implies that \begin{equation*}
m(t,x)=m_0(x)\exp^{-\int_0^tu_x(\tau,\,x)\,d\tau}.
\end{equation*}
By (\ref{eq3.7}), we see that for every $(t,x)\in[0,T)\times
\mathbb{R}$, $m(t,x)$ has the same sign as $m_0(x)$ does. This
completes the proof of the lemma. \end{proof}

\begin{Lemma}\label{lem3.4} Let $u_{0} \in H^{s}(\mathbb{R}),\,s>\frac{3}{2}$ be given.
If $m_{0}:=(u_{0}-u_{0,xx})\in L^{1}(\mathbb{R})$, then, as long
as the solution $u(t,\cdot)$ to Eq.(\ref{c0}) with initial data $u_{0}$
given by Theorem \ref{th3.1} exists, we have
$$
\int_{\mathbb{R}}u(t,x)dx=\int_{\mathbb{R}}u_{0}dx=
\int_{\mathbb{R}}m_{0}dx=\int_{\mathbb{R}}m(t,x)dx.
$$
\end{Lemma}
\begin{proof} Again it suffices to consider the case $s = 3$. Let $ T $ be the maximal time of existence of the solution $u$
to Eq.(\ref{c0}) with initial data $u_{0}\in H^{3}(\mathbb{R}).$

Note that $u_{0}=Q\ast m_{0}$ and $m_{0}=(u_{0}-u_{0,xx})\in
L^{1}(\mathbb{R})$. By Young's inequality, we get
\begin{equation*}
\| u_{0}\|_{L^{1}(\mathbb{R})}= \| Q\ast
m_{0}\|_{L^{1}(\mathbb{R})} \leq \| Q\|_{L^{1}(\mathbb{R})}\|
m_{0}\|_{L^{1}(\mathbb{R})}\leq \| m_{0}\|_{L^{1}(\mathbb{R})}.
\end{equation*}
Integrating Eq.(3.2) by parts, we get
$$
\frac{d}{dt}\int_{\mathbb{R}}u
dx=\int_{\mathbb{R}}\partial_{x}Q\ast
(\frac{1}{2}u^{2}_{x}-\frac{1}{2}u^{2})dx=0.
$$
It then follows that
$$
\int_{\mathbb{R}}u\,dx=\int_{\mathbb{R}}u_{0}\,dx.
$$
 Due to $m=u-u_{xx}$, we have
\begin{multline*}
\int_{\mathbb{R}}m\,dx=\int_{\mathbb{R}}u\,dx-
\int_{\mathbb{R}}u_{xx}dx=\int_{\mathbb{R}}u\,dx
\\=\int_{\mathbb{R}}u_{0}\,dx=\int_{\mathbb{R}}u_{0}\,dx-
\int_{\mathbb{R}}u_{0,xx}dx=\int_{\mathbb{R}}m_{0}\,dx.
\end{multline*}
This completes the proof of the lemma.
\end{proof}
We now present the first global existence result.
\begin{Theorem}\label{th3.4} Let $u_{0} \in H^{s}(\mathbb{R})\;s> \frac{3}{2}$
be given. If $m_0:=u_0-\partial_x^2 u_0\in L^{1}(\mathbb{R})$ is
nonnegative, then the corresponding solution to Eq.(\ref{c0+}) is
defined globally in time. Moreover, $I(u)=\int_{\mathbb{R}}u\,dx$
is a conservation law, and that for all $(t,x)\in \mathbb{R}_{+}\times \mathbb{R}$, we have\\
(i) $m(t,x)\geq 0$, $u(t,x)\geq 0$ and $$\|
m_{0}\|_{L^{1}(\mathbb{R})}=\| m(t)\|_{L^{1}(\mathbb{R})}=\|
u(t,\cdot) \|_{L^{1}(\mathbb{R})}=\|
u_{0}\|_{L^{1}(\mathbb{R})}.$$ (ii) $\|
u_{x}(t,\cdot)\|_{L^{\infty}(\mathbb{R})}\leq \|
u_0\|_{L^1(\mathbb{R})}$ and $$\|
u(t,\cdot)\|_{L^{\infty}(\mathbb{R})}\leq \frac{1}{\sqrt{2}}\| u(t,\cdot)\|_{1}\leq
\frac{\sqrt{2}}{2}e^{\frac{t}{2}\| u_{0}\|_{L^{1}(\mathbb{R})}} \|
u_{0}\|_{1}.$$
\end{Theorem}
\par
\begin{proof}  As we mentioned before that we only need to prove the above
theorem for $s=3$.  Let $T>0$ be the maximal existence time of the
solution $u$ with initial data $u_0\in H^{3}(\mathbb{R})$.

If $m_0(x)\geq 0$, then Lemma \ref{lem3.3} ensures that $m(t,x)\geq
0$ for all $t\in [0,T)$. Noticing that $u=Q\ast m$ and the
positivity of $Q$, we infer that $u(t,x)\geq 0$ for all $t\in
[0,T)$. By Lemma \ref{lem3.4}, we obtain
\begin{multline}\label{eq4.1}
-u_{x}(t,x)+\int_{-\infty}^{x}u(t,x)
dx=\int_{-\infty}^{x}(u-u_{xx})dx\\=\int_{-\infty}^{x}m dx
\leq\int_{-\infty}^{\infty}mdx=\int_{\mathbb{R}}m_{0}dx=
\int_{\mathbb{R}}u_{0}dx.
\end{multline}
Therefore, from (\ref{eq4.1}) we find that
\begin{equation}\label{eq4.2}
u_{x}(t,x)\geq -\int_{\mathbb{R}}u_{0}dx=-\|
u_{0}\|_{L^{1}(\mathbb{R})}, \quad \forall
(t,x)\in[0,T)\times\mathbb{R}.
\end{equation}
On the other hand, by $m(t,x)\geq 0$ for all $t\in [0,T)$, we obtain
\begin{equation*}
u_{x}(t,x)-\int_{-\infty}^{x}u\,
dx=-\int_{-\infty}^{x}(u-u_{xx})\,dx=-\int_{-\infty}^{x}m\,dx\leq
0.
\end{equation*}
By the above inequality and $u(t,x)\geq 0$ for all $t\in [0,T)$,
we get
\begin{equation}\label{eq4.3}
u_{x}(t,x)\leq\int_{-\infty}^{x}u\, dx\leq
\int_{\mathbb{R}}u\,dx=\int_{\mathbb{R}}u_{0}\,dx=\|
u_{0}\|_{L^{1}(\mathbb{R})}.
\end{equation}
Thus, (\ref{eq4.2}) and (\ref{eq4.3}) imply that
\begin{equation}
|u_{x}(t,x)|\leq
\|u_x(t,\cdot)\|_{L^{\infty}(\mathbb{R})}\leq \|
u_{0}\|_{L^{1}(\mathbb{R})}\quad \forall
(t,x)\in[0,T)\times\mathbb{R}.\end{equation} By Theorem \ref{th3.2} and
the above inequality, we deduce that $T=\infty$. Recalling finally
Lemma \ref{lem3.4}, we get assertion (i).
\par
Multiplying (\ref{c0}) by $u$ and integrating by parts, we obtain
\begin{equation}
\begin{split}
&\frac{1}{2}\frac{d}{dt}\int_{\mathbb{R}}\left(u^2(t,x)+u_x^2(t,x)\right)dx
=\int_{\mathbb{R}}(uu_xu_{xx}-u^2u_x)dx\\&=-\frac{1}{2}\int_{\mathbb{R}}u_x^3dx\leq
\frac{1}{2}\|u_x(t,\cdot)\|_{L^{\infty}(\mathbb{R})}\int_{\mathbb{R}}u_x^2dx.
\end{split}
\end{equation}
An application of Gronwall's inequality leads to
\begin{equation}
\int_{\mathbb{R}}\left(u^2(t,x)+u_x^2(t,x)\right)dx \leq
e^{t\|u_0\|_{L^{1}(\mathbb{R})}}\int_{\mathbb{R}}\left(u_0^2+u_{0,x}^2\right)dx.
\end{equation}
Consequently,
\begin{equation}\label{3.28}
\| u(t,\cdot)\|_{1}\leq e^{\frac{t}{2}\|
u_{0}\|_{L^{1}(\mathbb{R})}} \| u_{0}\|_{1}.
\end{equation}
On the other hand,
\begin{equation}\label{3.29}
u^2(t,x)= \int^{x}_{-\infty}uu_x\,dx-\int_{x}^{\infty}uu_x\,dx\leq
\frac{1}{2}\int_{\mathbb{R}}(u^2+u_x^2)dx=\frac{1}{2}\|
u(t,\cdot)\|^2_{1}.
\end{equation}
Combining (\ref{3.28}) with (\ref{3.29}), we obtain assertion (ii). This
completes the proof of the theorem.
\end{proof}
In a similar way to the proof of Theorem \ref{th3.4}, we can get the
following global existence result.
\begin{Theorem} \label{th3.5}Let $u_{0} \in H^{s}(\mathbb{R})\;s> \frac{3}{2}$ be given.
If $m_0:=u_0-\partial_x^2 u_0\in L^{1}(\mathbb{R})$ is
non-positive, then the corresponding solution to Eq.(\ref{c0}) is
defined globally in time. Moreover, $I(u)=\int_{\mathbb{R}}u\,dx$
is invariant in time, and that for all $(t,x)\in \mathbb{R}_{+}\times \mathbb{R}$, we have\\
(i) $m(t,x)\leq 0$, $u(t,x)\leq 0$ and $$\|
m_{0}\|_{L^{1}(\mathbb{R})}=\| m(t)\|_{L^{1}(\mathbb{R})}=\|
u(t,\cdot) \|_{L^{1}(\mathbb{R})}=\|
u_{0}\|_{L^{1}(\mathbb{R})}.$$
(ii) $\|
u_{x}(t,\cdot)\|_{L^{\infty}(\mathbb{R})}\leq \|
u_0\|_{L^1(\mathbb{R})}$ and $$\|
u(t,\cdot)\|_{L^{\infty}(\mathbb{R})}\leq \frac{\sqrt{2}}{2} \| u(t,\cdot)\|_{1}\leq
\frac{\sqrt{2}}{2}e^{\frac{t}{2}\| u_{0}\|_{L^{1}(\mathbb{R})}} \|
u_{0}\|_{1}.$$
\end{Theorem}

\section{Global weak solutions}

\noindent In this section, we present some results on global weak
solutions to characterize peakon solutions to (\ref{eq2.1}) for any
$\theta\in \mathbb{R}$ provided initial data satisfy certain sign
conditions.
\par
Let us first introduce some notations to be used in the sequel. We let
$M(\mathbb{R})$ denote the space of Radon measures on $\mathbb{R}$
with bounded total variation. The cone of positive measures is
denoted by $M^{+}(\mathbb{R})$. Let $BV(\mathbb{R})$ stand for the
space of functions with bounded variation and write
$\mathbb{V}(f)$ for the total variation of $f\in BV(\mathbb{R})$.
Finally, let $\{\rho_{n}\}_{n\geq 1}$ denote the mollifiers
$$
\rho_{n}(x):=\left(\int_{\mathbb{R}}\rho(\xi)d\xi\right)^{-1}n\rho(nx),\quad
x\in\mathbb{R},\;n\geq 1,
$$
where $\rho\in C_{c}^{\infty}(\mathbb{R})$ is defined by
$$
\rho(x):=\begin{cases}e^{\frac{1}{x^{2}-1}},\quad &\text{
for}\;\mid x\mid<1,\\0,\quad &\text{ for}\;\mid x\mid\geq 1.
\end{cases}
$$
\par
  Note that the $b-$equation for any $b\in\mathbb{R}$ has peakon solutions
  with corners at their peaks, cf. \cite{D-H-H1,H-S,H-S2}. Thus, the $\theta-$equation for any $\theta\in \mathbb{R}$ has also
  peakon solutions, see Example 6.1 below.
  Obviously, such solutions are not strong solutions to (\ref{eq2.1}) for any $\theta\in\mathbb{R}$. In order to provide a mathematical
  framework for the study of peakon solutions, we shall first give
  the notion of weak solutions to (\ref{eq2.1}).
  \par
 Equation (\ref{eq2.1}) can be written as
 $$
 u_t+ \theta u u_x +\partial_x(1-\partial_x^2)^{-1}B(u, u_x)=0, \quad B=(1-\theta)\frac{u^2}{2}+(4\theta-1)
 \frac{u_x^2}{2}.
 $$
If we set
  $$F(u):=\frac{\theta u^2}{2}+Q \ast \left[ (1-\theta)\frac{u^2}{2}+(4\theta-1) \frac{u_x^2}{2}\right],$$
  then the above equation takes the conservative form
\begin{equation}\label{fu}
u_{t}+F(u)_{x}=0,\quad u(0,x)=u_{0},\quad t>0,\;x\in\mathbb{R}.
\end{equation}
In order to introduce the notion of weak solutions to (\ref{fu}),
let $\psi\in C_{0}^{\infty}([0,T)\times \mathbb{R})$ denote the set
of all the restrictions to $[0,T)\times \mathbb{R}$ of smooth
functions on $\mathbb{R}^{2}$ with compact support contained in
$(-T,T)\times \mathbb{R}$.
\begin{Definition}
Let $u_{0}\in H^{1}(\mathbb{R})$. If $u$ belongs to
$L^{\infty}_{loc}([0,T);H^{1}(\mathbb{R}))$ and satisfies the
following identity
$$
\int_{0}^{T}\int_{\mathbb{R}}(u\psi_{t}+F(u)\psi_{x})dxdt+\int_{\mathbb{R}}u_{0}(x)\psi(0,x)dx=0
$$
for all $\psi\in C_{0}^{\infty}([0,T)\times \mathbb{R})$, then $u$
is called a weak solution to (\ref{fu}). If $u$ is a weak solution
on $[0,T)$ for every $T>0$, then it is called a global weak solution
to (\ref{fu}).
\end{Definition}
The following proposition is standard.
\begin{Proposition}
(i) Every strong solution is a weak solution. \\
(ii) If $u$ is a weak solution and $u\in
C([0,T);H^{s}(\mathbb{R}))\bigcap
C^{1}([0,T);H^{s-1}(\mathbb{R}))$, $s>\frac{3}{2}$, then it is a
strong solution.
\end{Proposition}

\par Referring to an approximation procedure used first for the
solutions to the Camassa-Holm equation \cite{C-M}, a partial
integration result in Bochner spaces \cite{M-N-R-R} and Helly's
theorem \cite{N} together with the obtained global existence results
and two useful a priori estimates for strong solutions, e.g.,
Theorem \ref{th1.1} and Theorems \ref{th3.4}-\ref{th3.5}, we may
obtain the following uniqueness and existence results for the global
weak solution to (\ref{fu}) for any $\theta\in\mathbb{R}$ provided
the initial data satisfy certain sign conditions.
\begin{Theorem}\label{th6.1}
Let $u_{0}\in H^{1}(\mathbb{R})$ be given. Assume that
$$(u_{0}-u_{0,xx})\in
M^{+}(\mathbb{R}).$$ Then (\ref{fu}) for any $\theta\in \mathbb{R}$
has a unique weak solution
$$u\in W^{1,\infty}(\mathbb{R}_{+}\times\mathbb{R})\cap
L^{\infty}_{loc}(\mathbb{R}_{+};H^{1}(\mathbb{R}))$$ with initial
data $u(0)=u_{0}$ and
$$(u(t,\cdot)-u_{xx}(t,\cdot))\in M^{+}(\mathbb{R})$$ is uniformly
bounded for all $t\in\mathbb{R}_+$.
\end{Theorem}

In the following, we only present main steps of the proof of the
theorem, a refined analysis could be done following those given in
\cite{E-Y} for the $b-$equation. \\

\noindent {\bf A sketch of existence proof of weak solutions}

Step 1. Given $u_{0}\in H^{1}(\mathbb{R})$ and $m_{0}:=u_{0}-
u_{0,xx}\in M^{+}(\mathbb{R})$. Then one can show that
$$
\|u_0\|_{L^1(\mathbb{R})}\leq \|m_0\|_{ M^{+}(\mathbb{R})}.
$$
Let us define $u_{0}^{n}:=\rho_{n}\ast u_{0}\in
H^{\infty}(\mathbb{R})$ for $n\geq 1$. Note that for all $n\geq 1$,
$$
m_{0}^{n}:=u^{n}_{0} -u^{n}_{0,xx} =\rho_{n}\ast (m_{0})\geq 0.
$$By Theorem \ref{th1.1} and
Theorem \ref{th3.4}, we obtain that there exists a unique strong
solution to (\ref{fu}),
$$u^{n}=u^{n}(.,u_{0}^{n})\in C([0,\infty);H^{s}(\mathbb{R}))\cap
C^{1}([0,\infty);H^{r-1}(\mathbb{R})),\quad\forall s\geq 3.
$$

Step 2. By a priori estimates in Theorem \ref{th1.1} and Theorem
\ref{th3.4}, Young's inequality and energy estimate for (\ref{fu}),
we may get
\begin{equation}\label{6.1}
\int_{0}^{T}\int_{\mathbb{R}}(
[u^{n}(t,x)]^{2}+[u^{n}_{x}(t,x)]^{2}+[u^{n}_{t}(t,x)]^{2})dxdt\leq
M,
\end{equation}
where $M$ is a positive constant depending only on $\theta$, $T$,
$\|Q_x\|_{L^2(\mathbb{R})}$, and $ \| u_{0}\|_{1}$. It then follows
from (\ref{6.1}) that the sequence $\{u^{n}\}_{n\geq 1}$ is
uniformly bounded in the space $H^{1}((0,T)\times\mathbb{R})$. Thus,
we can extract a subsequence such that
\begin{equation}
u^{n_{k}}\rightharpoonup u \quad \text{weakly
in}\;H^{1}((0,T)\times\mathbb{R})\;\,\text{for}\;n_{k}\rightarrow\infty
\end{equation}
and
\begin{equation}
u^{n_{k}}\longrightarrow u \quad \text{a.e.
on}\;(0,T)\times\mathbb{R}\;\,\text{for}\;\,n_{k}\rightarrow\infty,
\end{equation}
for some $u\in H^{1}((0,T)\times\mathbb{R})$. From Theorem
\ref{th1.1} (i) and the fact $\|u_0^n\|_1\leq \|u_0\|_1$ we see that
for any fixed $t \in (0, T)$, the sequence $u_x^{n_k}(t, \cdot)\in
BV(\mathbb{R})$ satisfies
$$
\mathbb{V}[u_x^{n_k}(t, \cdot)]\leq 2\|m_0\|_{M(\mathbb{R})}.
$$

Step 3. This, when applying Helly's theorem, cf. \cite{N}, enables
us to conclude that there exists a subsequence, denoted still by
$\{u_{x}^{n_{k}}(t,\cdot)\}$, which converges to the function
$u_{x}(t,\cdot)$ for a.e. $t\in(0,T)$.  A key energy estimate is of
the form
$$
\|B(u^n, u^n_x)\|_{L^2(\mathbb{R})}\leq C(\|u_0\|_1),
$$
which ensures $B$ admits a weak limit. This when combined with the
fact that $(u^n, u_x^n)$ converges to $(u, u_x)$ as well as $Q_x\in
L^2$ leads to the assertion that $u$ satisfies (\ref{fu}) in
distributional sense.

Step 4. From equation (\ref{fu}) we see that $u^{n_k}_t(t, \cdot)$
is uniformly bounded in $L^2(\mathbb{R})$, and $\|u^{n_k}(t,
\cdot)\|_1$ is uniformly bounded for all $t\in (0, T)$.  This
implies that the map $t |\to u^{n_k}_t(t, \cdot) \in
H^1(\mathbb{R})$ is weakly equi-continuous on $[0, T]$. Recalling
the Arzela-Ascoli theorem and a priori estimates in Theorem
\ref{th1.1} and Theorem \ref{th3.4}, we may prove
$$u\in L_{loc}^{\infty}(\mathbb{R}_{+}\times \mathbb{R})\cap
L_{loc}^{\infty}(\mathbb{R}_{+};H^1(\mathbb{R}))\quad \text{and}
\quad u_{x}\in L^{\infty}(\mathbb{R}_{+}\times \mathbb{R}).$$

Step 5. Since $u$ solves (\ref{fu}) in distributional sense, we
have
\begin{equation*}
\rho_{n}\ast u_{t}+\rho_{n}\ast(\theta
uu_{x})+\rho_{n}\ast\partial_{x}Q\ast
((1-\theta)\frac{u^{2}}{2}+(4\theta-1)\frac{u_x^{2}}{2})=0,
\end{equation*}
for a.e. $t\in\mathbb{R}_{+}$. Integrating the above equation with
respect to $x$ on $\mathbb{R}$ and then integrating by parts, we
obtain
\begin{equation*}
\frac{d}{dt}\int_{\mathbb{R}}\rho_{n}\ast u\,dx=0.
\end{equation*}
By a partial integration result in Bochner spaces \cite{M-N-R-R}
and Young's inequality, we may prove that
$$
\int_{\mathbb{R}}u(t,\cdot)dx=\lim_{n\rightarrow\infty}\int_{\mathbb{R}}\rho_{n}\ast
u(t,\cdot)dx=\lim_{n\rightarrow\infty}\int_{\mathbb{R}}\rho_{n}\ast
u_{0}dx=\int_{\mathbb{R}}u_{0}dx.
$$
Using the above conservation law, we get
\begin{equation*}
\begin{split}
\parallel u(t,\cdot)-u_{xx}(t,\cdot)\parallel_{M(\mathbb{R})}&\leq
\parallel u(t,\cdot)\parallel_{L^{1}(\mathbb{R})}+\parallel u_{xx}(t,\cdot)\parallel_{M(\mathbb{R})}
\\&\leq \parallel
u_{0}\parallel_{L^{1}(\mathbb{R})}+2\parallel
m_{0}\parallel_{M(\mathbb{R})}\leq 3\parallel
m_{0}\parallel_{M(\mathbb{R})}, \end{split}\end{equation*} for
a.e. $t\in\mathbb{R}_{+}$. Note that
$u^{n_{k}}(t,x)-u^{n_{k}}_{xx}(t,x)\geq 0$ for all
$(t,x)\in\mathbb{R}_{+}\times\mathbb{R}$. Then the above
inequality implies that $(u(t,\cdot)-u_{xx}(t,\cdot))\in
M^{+}(\mathbb{R})$ for a.e. $t\in \mathbb{R}_{+}$.
\par
Since $u(t,x)=Q\ast (u(t,x)-u_{xx}(t,x))$, it follows that
\begin{equation*}
\begin{split}
\mid u(t,x)\mid&=\mid Q\ast (u(t,x)-u_{xx}(t,x))\mid\\&\leq
\parallel Q\parallel_{L^{\infty}(\mathbb{R})}\parallel
u(t,\cdot)-u_{xx}(t,\cdot)\parallel_{M(\mathbb{R})}\leq
\frac{3}{2}\parallel m_{0}\parallel_{M(\mathbb{R})}.
\end{split}\end{equation*} This shows that $u(t,x)\in
W^{1,\infty}(\mathbb{R}_{+}\times\mathbb{R})$ in view of Step 4.
This proves the existence of global weak solutions to (\ref{fu}).
\bigskip\par
\noindent {\bf Uniqueness of the weak solution}\\
Let $$u,\,v\in W^{1,\infty}(\mathbb{R}_{+}\times\mathbb{R})\cap
L_{loc}^{\infty}(\mathbb{R}_{+};H^{1}(\mathbb{R}))$$ be two global
weak solutions of (\ref{fu}) with initial data $u_0$. Set
$$
N:=\sup_{t\in\mathbb{R}_{+}}\{\| u(t,\cdot)-u_{xx}(t,\cdot)\|_{
M(\mathbb{R})}+\| v(t,\cdot)-v_{xx}(t,\cdot)\|_{ M(\mathbb{R})}\}.
$$
From Step 5, we know that $N<\infty$. Let us set
$$
w(t,\cdot)=u(t,\cdot)-v(t,\cdot),\quad
(t,x)\in\mathbb{R}_{+}\times\mathbb{R},
$$
and fix $T>0$. Convoluting Eq.(\ref{fu}) for $u$ and $v$ with
$\rho_{n}$ and with $\rho_{n,x}$ respectively, using Young's
inequality and following the procedure described on page 56-59 in
\cite{C-M}, we may deduce that
\begin{equation}\label{6.5}
\frac{d}{dt}\int_{\mathbb{R}}\mid \rho_{n}\ast
w\mid\,dx=C\int_{\mathbb{R}}\mid\rho_{n}\ast
w\mid\,dx+C\int_{\mathbb{R}}\mid\rho_{n}\ast
w_{x}\mid\,dx+R_{n}(t),
\end{equation}
and
\begin{equation}\label{6.6}
\frac{d}{dt}\int_{\mathbb{R}}\mid \rho_{n}\ast
w_{x}\mid\,dx=C\int_{\mathbb{R}}\mid\rho_{n}\ast
w\mid\,dx+C\int_{\mathbb{R}}\mid\rho_{n}\ast
w_{x}\mid\,dx+R_{n}(t),
\end{equation}
for a.e. $t\in [0,T]$ and all $n\geq 1$, where $C$ is a generic
constant depending on $\theta$ and $N$, and that $R_n(t)$
satisfies
\begin{equation*}
\begin{cases} \lim\limits_{n\rightarrow\infty}R_{n}(t)=0\\ \mid
R_{n}(t)\mid\leq K(T),\quad n\geq 1,\quad t\in[0,T].
\end{cases}
\end{equation*}
Here $K(T)$ is a positive constant depending on $\theta$, $T$, $N$
and the $H^{1}(\mathbb{R})$-norms of $u(0)$ and $v(0)$.
\par
Summing (\ref{6.5}) and (\ref{6.6}) and then using Gronwall's
inequality, we infer that
\begin{multline*}
\int_{\mathbb{R}}\bigl(\mid \rho_{n}\ast w\mid+\mid \rho_{n}\ast
w_{x}\mid\bigr)(t,x)\,dx\leq \int_{0}^{t}e^{2C(t-s)}R_{n}(s)ds+\\
e^{2Ct}\int_{\mathbb{R}}\bigl(\mid \rho_{n}\ast w\mid+\mid
\rho_{n}\ast w_{x}\mid\bigr)(0,x)\,dx,
\end{multline*}
for all $t\in[0,T]$ and $n\geq 1$. Note that $w=u-v\in
W^{1,1}(\mathbb{R})$. Using Lebesgue's dominated convergence
theorem, we may deduce that for all $t\in [0,T]$
$$
\int_{\mathbb{R}}\bigl(\mid w\mid+\mid
w_{x}\mid\bigr)(t,x)\,dx\leq e^{2Ct}\int_{\mathbb{R}}\bigl(\mid
w\mid+\mid w_{x}\mid\bigr)(0,x)\,dx.
$$
Since $w(0)=w_{x}(0)=0$, it follows from the above inequality that
$u(t,x)=v(t,x)$ for all $(t,x)\in [0,T]\times\mathbb{R}$. This
proves the uniqueness of the global weak solution to (\ref{fu}).

\bigskip
\par
In a similar way to the proof of Theorem 6.1, we can get the
following result.

\begin{Theorem}\label{th6.2}
Let $u_{0}\in H^{1}(\mathbb{R})$ be given. Assume that
$$(u_{0,xx}-u_{0})\in
M^{+}(\mathbb{R}).$$ Then (\ref{fu}) for any $\theta\in \mathbb{R}$
has a unique weak solution
$$u\in W^{1,\infty}(\mathbb{R}_{+}\times\mathbb{R})\cap
L^{\infty}_{loc}(\mathbb{R}_{+};H^{1}(\mathbb{R}))$$ with initial
data $u(0)=u_{0}$ and
$$(u_{xx}(t,\cdot)-u(t,\cdot))\in M^{+}(\mathbb{R})$$ is uniformly
bounded for all $t\in\mathbb{R}_+$.
\end{Theorem}
\begin{Remark} Theorems \ref{th6.1}-\ref{th6.2} cover the recent results for global weak solutions of the
Camassa-Holm equation in \cite{C-M} and the Degasperis-Procesi
equation in \cite{Y4}.
\end{Remark}
\begin{Example} (Peakon solutions) Consider (\ref{eq2.1}) for any $\theta\in\mathbb{R}$.
Given the initial datum $u_{0}(x)=ce^{-\mid x\mid}$,
$c\in\mathbb{R}$. A straightforward computation shows that
$$u_{0}-u_{0,xx}=2c\,\, \delta(x)\in
M_{+}(\mathbb{R})\quad \mbox{if}\quad c\geq 0$$ and
$$u_{0,xx}-u_{0}=-2c\,\, \delta(x)\in
M_{+}(\mathbb{R})\quad \mbox{if}\quad c< 0.$$ One can also check
that
$$u(t,x)=ce^{-\mid
x-ct\mid}$$ satisfies (\ref{eq2.1}) for any $\theta\in\mathbb{R}$
in distributional sense. Theorems \ref{th6.1}-\ref{th6.2} show
that $u(t,x)$ is the unique global weak solution to (\ref{eq2.1})
for any $\theta\in\mathbb{R}$ with the initial data $u_{0}(x)$.
This weak solution is a peaked solitary wave which is analogues to
that of the $b-$equation, cf. \cite{E-Y}.
\end{Example}

\noindent\textbf{Acknowledgments.} The authors gratefully
acknowledge the support of the basic research program ``Nonlinear
Partial Differential Equations" at the Center for Advanced Study
at the Norwegian Academy of Science and Letters, where this work
was performed during their visit in December of 2008. Liu's
research was partially supported by the National Science
Foundation under the Kinetic FRG Grant DMS07-57227 and the Grant DMS09-07963.
Yin's research was partially supported by NNSFC (No. 10971235 and No. 10531040), RFDP
(No. 200805580014) and NCET-08-0579.

\end{document}